# Social cognitive optimization with tent map for combined heat and power economic dispatch

Jiaze Sun [a,b] *, Yang Li [c]

[a] *School of Computer Science and Technology, Xi'an University of Posts and Telecommunications, Xi'an 710121, China*

[b] *Shaanxi Key Laboratory of Network Data Intelligent Processing, Xi'an University of Posts and Telecommunications, Xi'an 710121, China*

[c] *School of Electrical Engineering, Northeast Electric Power University, Jilin 132012, China*

**ABSTRACT:** Combined Heat and Power Economic Dispatch (CHPED) problem is a sophisticated constrained nonlinear optimization problem in a heat and power production system for assigning heat and power production to minimize the production costs. To address this challenging problem, a novel Social Cognitive Optimization algorithm with Tent map (TSCO) is presented for solving the CHPED problem. To handle the equality constraints in heat and power balance constraints, Adaptive Constraints Relaxing (ACR) rule is adopted in constraint processing. The novelty of our work lies in the introduction of a new powerful TSCO algorithm to solve the CHPED issue. The effectiveness and superiority of the presented algorithm is validated by carrying out two typical CHPED cases. The numerical results show that the proposed approach has better convergence speed and solution quality than all other existing state-of-the-art algorithms.

## 1. Introduction

Despite the reducing fossil fuels and the increasing need for electric energy, the efficiency of the conversion of fossil fuels into electricity is still less than 60%. Plenty of thermal energy as wasted heat is directly released into the environment, which leads to severe environmental pollution problems. To make full use of wasted heat, Combined Heat and Power (CHP) units are introduced to simultaneously produce electricity and heat energy in the power production by conventional thermal plants. CHP units can generate conversion efficiency above than 90% through using wasted heat in thermal plant [1], so the cogeneration units, power units, and heat units are combined to satisfy power and heat demands in CHP plants.

To achieve the comprehensive utilization of CHP units, Combined Heat and Power Economic Dispatch (CHPED) optimization is implemented for obtaining an optimal power and heat generation scheduling, which minimizes the overall cost of supplying power and heat demand while satisfying system all inequality and equality constraints. The CHPED optimization problem is a multimodal, non-convex and constrained nonlinear programming problems [2]. Considering coupling of power and heat, the valve point effects and transmission losses between multiple plants, the sophistication of the CHPED problem increases further [3].

The literatures [1, 3, 4] show that various optimization methods have been introduced to solve the CHPED problem. Several traditional optimization approaches, such as gradient descent approach, nonlinear optimization method, have been introduced to work out the CHPED problems, but the traditional deterministic approaches are very difficult to precisely obtain the optimal solution because the model of the sophisticated CHPED problem is always non-smooth with many complex equality and inequality constraints. Moreover, this constrained nonlinear global optimization problem has proven to be non-deterministic polynomial-time hard (NP-hard) [1], up to now there is still no effective deterministic method to solve them in practice.

When solving the multimodal, non-convex and non-differentiable optimization problem, Swarm Intelligence (SI) optimization algorithms have been identified to have effective performance. Recent research demonstrates that the NP-hard problems can be approximately settled by metaheuristic or heuristic algorithms [4, 6, 7]. Many effective meta-heuristic or heuristic algorithms are implemented to search the optimal solution for the CHPED optimization problem [1]. The metaheuristic algorithms applied in CHPED optimization problem from the literatures include

Genetic Algorithm (GA) [8-9], Particle Swarm Optimization (PSO) algorithm [10-11], Ant Colony Optimization (ACO)algorithm [12], Bee Colony Optimization (BCO) [13], Cuckoo Search (CS)algorithm [14], Grey Wolf Optimization (GWO) algorithm [15], Artificial Immune System (AIS) algorithm [16], Firefly Algorithm (FA) [17], Harmony Search(HS) [18] , Differential Evolution (DE) [19] , Fish School Search (FSS) [20], Invasive Weed Optimization(IWO) algorithm [21], Group Search Optimization(GSO) [22] and Teaching Learning Based Optimization (TLBO) [23]. Although the above swarm intelligence optimization algorithms have had acceptable results to solve the CHPED optimization problem, the global optimization solution and its global convergence are usually difficult to be guaranteed. Moreover, the other notable meta-heuristic algorithms have not been explored.

It is known that human has better social intelligence and higher fitness than the other swarm. Vitalized by the human cognitive process, Social Cognitive Optimization (SCO) algorithm [24] was devised as a heuristic swarm algorithm. SCO algorithm includes many learning agents, which have symbolizing ability and vicarious ability. Since SCO algorithm makes optimum use of the overall social knowledge, it has better optimal ability than many other well-known optimization algorithms in many fields [25], such as GA，ACO and PSO.

To obtain a powerful, stable and global SCO algorithm for the CHPED problem, we combine the social cognitive optimization algorithm with Tent map and Adaptive Constraints Relaxing (ACR) for addressing the challenging CHPED problem. The main contributions of the paper are as follows:

(1) A novel social cognitive optimization algorithm based on Tent map and Adaptive Constraints Relaxing is proposed for the CHPED problem in the heat and power production system.

(2) Integrating Tent map into SCO initialization and neighborhood searching improves uniform distribution and increases the ergodicity of the SCO with less computation. Furthermore, ACR rule is adopted in constraint processing to effectively handle the equation constraints in heat and power balance constraints.

(3) The numerical experimental results that assess the effectiveness and convergence of the proposed novel social cognitive optimization with two benchmarks testing projects are compared with the other classic heuristic algorithm for the CHPED problem.

The remainder of this paper is structured as follows. The CHPED formulation is detailed in a constraint non-linear program model in Section 2. The specific description of the novel SCO with tent map for the CHPED problem is provided in section 3. In section 4 two classic test cases are exercised to assess the effectiveness and superiority of the novel SCO. In Section 5 the conclusions are shown.

## 2. Optimization formulation of CHPED Problem

A CHP system includes power units, cogeneration units, and heat units. Mathematically speaking, the CHPED problem can be formulated as constrained Nonlinear Programming Problems (NLPs) aiming to minimize total cost of supplying power and heat demand.

### 2.1 cost function

In practice, when steam admission valve starts to open, an intense fuel loss raising will enhance the fuel cost because of the wire drawing outcomes, which results in a non-smooth and non-convex fuel cost function. If the valve-point effects for power units are neglected, the fuel cost function is closely a quadratic function, which would lead an inaccuracy to the CHPED problem. Consequently, the objective function of the practical CHPED problem with the consideration of the valve-point effect is a superposition of quadratic and sinusoidal functions, which would increase the non-smooth and non-convex characteristic [4]. The cost function of power unit considering the valve-point effect can be represented as follows:

$$\text{cost}_i(P_i^p) = a_i(P_i^p)^2 + b_i p_i^p + c_i + |d_i \sin(e_i(P_i^{P_{\min}} - P_i^p))| \qquad (1)$$

In which, $\text{cost}_i(P_i^p)$ is the cost of power unit $i$, $P_i^p$ is the power generated by the power unit $i$, $a_i, b_i$, and $c_i$ are the cost parameters of power unit $i$, $d_i$ and $e_i$ are the cost parameters of the valve-point effect. $P_i^{P_{\min}}$ is the lower power boundaries of the power unit $i$ in MW.

Additionally, the cost function of cogeneration unit is stated as follows:

$$\text{cost}_j(P_j^c, H_j^c) = a_j(P_j^c)^2 + b_j P_j^c + c_j + d_j(H_j^c)^2 + e_j H_j^c + f_j H_j^c P_j^c \qquad (2)$$

where $\text{cost}_j(P_j^c, H_j^c)$ is the cost function of cogeneration unit $j$, and $a_j, b_j, c_j, d_j, e_j$ and $f_j$ are the cost parameters of cogeneration unit $j$. $P_j^c$ is the power generated by the cogeneration unit $j$, $H_j^c$ is the heat generated by the cogeneration unit $j$.

Moreover, the double-dip function can model the cost of heat unit, which is stated as follows:

$$\text{cost}_k(H_k^h) = a_k(H_k^h)^2 + b_k H_k^h + c_k \qquad (3)$$

in which, $\text{cost}_k(H_k^h)$ is the cost function of heat unit $k$, and $H_k^h$ is the heat generated by heat unit $k$. $a_k, b_k$, and $c_k$ are the cost parameters of heat unit $k$.

## 2.2 objective function

According to the above cost functions, the objective function of the CHPED problem can be written as

$$\text{minmize} \quad f_{\text{cost}} = \sum_{i=1}^{N_p} \text{cost}_i(P_i^p) + \sum_{j=N_p+1}^{N_p+N_c} \text{cost}_j(P_j^c, H_j^c) + \sum_{k=N_p+N_c+1}^{N_p+N_c+N_h} \text{cost}_k(H_k^h) \qquad (4)$$

where $f_{\text{cost}}$ is the total cost (\$/h). $N_p, N_c$, and $N_h$ are the respective number of power units, cogeneration units, and heat units. The variables $i, j$ and $k$ are respectively the power unit number, cogeneration unit number, and heat unit number.

## 2.3 constraint

Since the power transmission loss exists factually in the power transmission system, network transmission loss is an important factor in the CHPED problem. In general, the network losses can be calculated through power generation of all units which is known as B-matrix approach. Considering transmission line losses using B-matrix method is reflected as follows:

$$P_{Loss} = \sum_{i=1}^{N_p}\sum_{j=1}^{N_p} P_i^p B_{ij} P_j^p + \sum_{i=1}^{N_p}\sum_{j=N_p+1}^{N_p+N_c} P_i^P B_{ij} P_j^c + \sum_{i=N_p+1}^{N_p+N_c}\sum_{j=N_p+1}^{N_p+N_c} P_i^c B_{ij} P_j^c + \sum_{i=1}^{N_p} B_{0i} P_i^p + \sum_{i=N_p+1}^{N_p+N_c} B_{0i} P_j^c + B_{00} \qquad (5)$$

in which, $P_{Loss}$ is the power transmission of the system. Matrix B is the coefficients of the transmission power loss with dimension $(N_p+N_c) \times (N_p+N_c)$. $B_0$ is a vector with dimension $(N_p+N_c)$, and $B_{00}$ is a real constant number.

The equality constraints representing the power and heat demands and inequality constraints representing the capacity boundaries are given as follows:

**Power balance**

$$\sum_{i=1}^{N_p} P_i^p + \sum_{j=N_p}^{N_p+N_c} P_j^c - (P_d + P_{Loss}) = 0 \qquad (6)$$

In which, $P_d$ and $P_{Loss}$ are the indexes of power demand and the power transmission losses.

**Heat balance**

$$\sum_{j=N_p}^{N_p+N_c} H_j^c + \sum_{k=N_p+N_c}^{N_p+N_c+N_h} H_k^h - H_d = 0 \qquad (7)$$

$H_d$ is the indicator of heat demand.

**Capacity boundary**

$$P_i^{P\min} \leq P_i^p \leq P_i^{P\max} \qquad i=1,...,N_p \qquad (8)$$

$$P_j^{c\min}(H_j^c) \leq P_j^c \leq P_j^{c\max}(H_j^c) \qquad j=N_p+1,...,N_p+N_c \qquad (9)$$

$$H_j^{c\min}(P_j^c) \leq H_j^c \leq H_j^{c\max}(P_j^c) \qquad j=N_p+1,...,N_p+N_c \qquad (10)$$

$$H_k^{h\min} \leq H_k^h \leq H_k^{h\max} \qquad k=N_p+N_c+1,...,N_p+N_c+N_h \qquad (11)$$

The upper and lower power boundaries of the power unit $i$ are $P_i^{P\max}$ and $P_i^{P\min}$ (MW). The upper and lower power boundaries of cogeneration unit $j$ are $P_j^{c\min}(H_j^c)$ and $P_j^{c\max}(H_j^c)$ (MW), which are functions of heat $H_j^c$. The upper and lower heat boundaries of cogeneration unit $j$ are $H_j^{c\min}(P_j^c)$ and $H_j^{c\max}(P_j^c)$ (MWth) which are functions of power in the cogeneration unit $j$. The upper and lower heat boundaries of the heating unit $k$ are $H_k^{h\max}$ and $H_k^{h\min}$ (MWth).

For ease of description, the CHPED problem as a typical constrained NLPs can be converted to the following traditional form:

$$\begin{aligned} min \quad & f(X) \\ s.t \quad & g(X) \leq 0 \\ & h(X) = 0 \end{aligned} \qquad (12)$$

where $\mathbf{X}=(x_1,...,x_q,...,x_n)$ $(1 \leq q \leq n, q \in \mathbb{Z})$ are the decision variables. $f(\mathbf{X})$ is the minimum objective function. $g(\mathbf{X})=[g_1(X),...,g_k(X)]^T$ is a vector of $k$ inequality constraints, and $h(\mathbf{X})=[h_1(X),...,h_m(X)]^T$ is a vector of m equality constraints.

3 **Social cognitive optimization algorithm with tent map**

3.1 **Social cognitive optimization (SCO)**

The social cognitive theory argues that the mankind studies by observing others with the surrounding, behavior, and cognition as the main interactive development factors. Human learning owns the capability to symbolize, learn from others, plan alternative strategies, regulate one's own behavior, and engage in self-reflection. Social cognitive experts argue that mankind has greater social intelligence and fitness than insect society. By adopting social cognitive theory in an artificial system, SCO algorithm was proposed in literature [24]. In SCO iteration process, knowledge library is composed of many knowledge points, which include the locations and fitness values. Learning agents take part in observational learning by the local neighbor searching. Referring to $x_1$, the local neighbor searching for $x_2$ is to obtain another point $x'$, which is formulated for $d$ dimension as follows.

$$x' = U(x_1, 2x_2 - x_1) \tag{13}$$

where $U(u,v)$ is a uniform distribution, which is always produced by the linear congruential approach in the interval $[u,v]$. Furthermore, it is necessary to ensure that $x'$ satisfy all the constraints. To select an appropriate point from the special set, tournament selection chooses the better point for neighborhood searching from the knowledge library and chooses the worse point for refreshing the knowledge library. More detailed procedure about the SCO algorithm can be found in literature [24].

### 3.2 Tent map

Chaos is a type of ubiquitous nonlinear phenomena in lots of actual systems. Chaotic movement can reach every state in certain scale according to its own regularity and ergodicity, which is better than a simple stochastic search algorithm. The chaotic search algorithm is featured with randomicity, ergodicity, and regularity. So, in many optimization algorithms, it is usually introduced to initialize the initial solution or to iterate local search to enhance the ergodicity of the solution and accelerate the global optimal convergence.

Chaos search algorithm has many implement models including Tent map, Kent map, and Logistic map. Tent map iteration is faster and more suitable for a computer than Logistic map and Kent map [25]. Meanwhile, the distribution of Tent map is very even, and its initial sensitivity is weak. After Bernoulli shift transformation, Tent map can be stated in the interval [0, 1]

$$x_{k+1} = \begin{cases} 2x_k, & 0 \le x_k \le 1/2 \\ 2(1-x_k), & 1/2 < x_k \le 1 \end{cases} \tag{14}$$

In Eq. (14), there are minor rotation points in tent map, such as four cycles: 1/5, 2/5, 3/5 and 4/5, and unstable cycle points, such as 0, 1/4, 1/2, 3/4 and 1, which will fall into the fixed point 0.

The chaotic series in the interval $(V_l, V_r)$ based on Tent map is produced as follows:

Step1: give an initial value which must avoid taking a few special values of 0, 1/4, 1/2, 3/4 and 1. Through $m$ (about 300) times iteration, numerical values $x_m$ in (0, 1) are gotten. In an iterative procedure, if $x_i = \{0, 0.25, 0.5, 0.75, 1\}$ or $x_i = x_{i-k}$  $k = \{1,2,3,4\}$, then $x_i = x_{i-1} + \varepsilon$ to avoid variable rotation points and minor rotation points. $\varepsilon$ is a very small real, for example $10^{-3}$.

Step2: calculating chaotic variables in each iteration according to Eq. (15)

$$f(x) = V_l + x_m \cdot |V_r - V_l| \tag{15}$$

where $V_r$ is the right value of the interval, and $V_l$ is the left value of the interval.

### 3.3 Coordinating Tent map and SCO

To obtain even distribution and improve the ergodicity of SCO initialization and local searching, chaotic search algorithm [26] is applied to library initialization and local searching of SCO instead of the common random approach in the classical SCO algorithm. In the social agents' initialization phrase, the initial position of agents is generated by the chaotic search algorithm in the feasible space rather than the stochastic method. In the local neighbor searching process, the local neighbor searching of the learning agents adopts the chaotic search algorithm to generate the new agent's position. Because the Tent map algorithm has more benefits than Kent map and Logistic map, the common random approach will be shifted by Tent map searching algorithm. Referring to $x_1$, the formulation of the local searching for $x_2$ is

$$x^{'} = \text{Tent}(x_1, 2x_2 - x_1) \qquad (16)$$

where $Tent(u,v)$ is a stochastic value generated by Tent map within $[u, v]$.

### 3.4 TSCO for CHPED problem

The decision variables in the CHPED problem are power and heat dispatch outputs values. The position of learning agent represents a feasible CHPED output value scheme in constraint space and its fitness function. The fitness includes two parts: CHPED objective function and constraints violation value. Each agent shows the possible solution of the CHPED problem, and all the agents form knowledge library. The *i*th knowledge point is represented as:

$$L_i = [P_{i1}^{\,p},...,P_{iN_p}^{\,p}, P_{i(Np+1)}^{\,c},...,P_{i(N_p+N_c)}^{\,c}, H_{i(Np+1)}^{\,c},...,H_{i(N_p+N_c)}^{\,c}, H_{i(N_p+N_c+1)}^{\,c},...,H_{i(N_p+N_c+N_h)}^{\,c}]$$

The length of the knowledge point is $N_p + N_c + N_h$. Learning agents in possession of knowledge points operate observational learning by the local Tent map searching and model selection via tournament selection. Fig. 1 illustrates the flowchart of Social Cognitive Optimization with Tent map (TSCO). The TSCO for CHPED problem is depicted as follows:

Step 1: initialization

1) Set parameters: the number of knowledge points $N_L$, the number of learning agents $N_a$, maximum iteration times $T$, the vicarious width $\tau_W$, the tournament width $\tau_B$.
2) Create $N_L$ knowledge points in knowledge library ($K$) with Tent map, and calculate their fitness values including objective function value and constraints violation value, then preserve the global optimal point $G_P$.
3) Every learning agent is assigned to a known point in $K$ at random, but not reduplicatively.

Step 2: every learning agent $L$:

1) Tournament selection: Choose the best-known point $T_P$ from random $\tau_B$ knowledge points in knowledge library $K$, not same with $L$ itself.
2) Observational learning: After the point $T_P$ is confronted with the point $L$, the better point is selected as a central point to generate the new point $T_O$ according to Eq. (16) with Tent map by referring to the worse. If the point $T_O$ is better than $G_P$, $T_O$ is assigned to $G_P$.
3) Library refreshment: Choose the worst point $T_w$ from random $\tau_W$ points in knowledge library $K$. If the fitness value of a point $T_O$ is better than the fitness value of $T_w$, the point $T_O$ will shift point $T_w$.

Step 3: iterate Step 2 until the stop condition is satisfied. The overall calculation time is $T_e = N_L + N_a * T$.

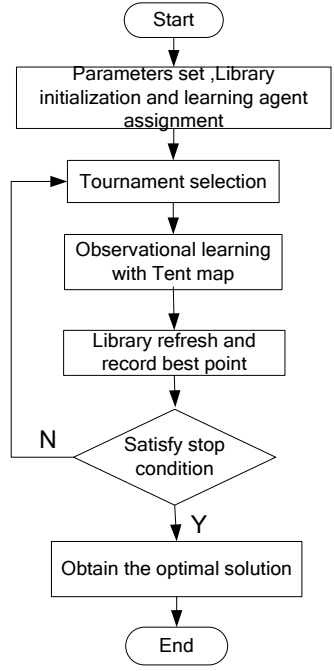

FIG. 1 Flowchart of TSCO

### 3.5 Constraint handling

Most literature applies a penalty function method to handle the constraints for the CHPED problem [1]. But the penalty handling strategy needs the tough penalty coefficient settings and feasible initial point. Basic Constraint Handing (BCH) strategy is introduced to handle the constraints in the swarm intelligence optimization algorithm. However, BCH strategy is difficult to cope with equality constraints in the swarm, since the violation values of the equality constraints are usually very small and easy to fall into the ridge function class landscape [27].

Because there are some power and heat balance equations constraints in CHPED, we adopt the relaxed quasi-feasible space to approximately meet the equation in the early iteration period. The actual feasible space is adaptively entered in the later iterations period [27]. The following formulation is introduced to calculate the constraint violation values:

$$F_{CON}(X) = \max\{\varepsilon_R, \sum_{i=1}^{m}\max\{0, h_i(X)\} + \sum_{j=1}^{k}\max\{0, g_j(X)\}\} \tag{17}$$

where $\varepsilon_R \geq 0$ is the relaxing threshold value. In the iterations, $\varepsilon_R$ constringes adaptively to zero based on the fundamental ratio-keeping latent rules and forcing latent rules [28].

### 4. Case studies

In this section, the performance of TSCO for resolving the CHPDE problems in the field of quality of solution and convergence speed is studied. Two typical test cases are chosen from literature [1] for comparison, which are used to evaluate the algorithm in most literature of studying the CHPDE problems. The first case is the classical simple representation of the CHPDE problem, and the second case is the classical sophisticated representation. The program has been implemented on the Eclipse Kepler SR1 IDE in Java language and executed on an Intel i7-5500U CPU@ 2.4 GHz 2.39 GHz PC with 8 GB RAM (thank Dr. Xie for sharing the source code of traditional SCO algorithm in the website: http://www.wiomax.com/team/xie/social-cognitive-optimization-sco-project-portal/). If different values are not clearly described in the test cases, the parameters of the TSCO for all test cases are set as follows: $N_L$ =150, $N_a$ =50, $T$ =100, $\tau_W$ =4, $\tau_B$ =2. Cost is in \$, heat output is in MWth, and power output is in MW in all the test cases. The numerical experimental results using traditional SCO and TSCO in this paper make a comparison with the typical literature [1] reports in the field of convergence speed and solution quality.

### 4.1 Test case 1

This typical test case includes four units which are one conventional power unit, one heat unit and two cogeneration units, and it is the classical representation of simple CHPDE problems and used to evaluate the algorithms in most literature of the CHPDE problems. In order to compare the proposed algorithm with the typical algorithms, we choose the classical simple test case as the benchmark, and we also do not considerate the valve effect as the same as in the other referred literatures. The test case power demand $P_d$ is 200MW and the test case heat demand $H_d$ is 115MWth. In this test case, the power transmission loss and valve-point effects are ignored as the same as in all the other literatures. The test system has six decision variables $(P_1, P_2, P_3, H_2, H_3, H_4)$. The fuel cost formulation of the power unit 1, co-generation unit 2, co-generation unit 3 and heat unit 4 are given:

$$\text{cost}_1(P_1) = 50P_1$$

$$\text{cost}_2(P_2, H_2) = 2650 + 14.5P_2 + 0.0345P_2^2 + 4.2H_2 + 0.03H_2^2 + 0.031P_2H_2$$

$$\text{cost}_3(P_3, H_3) = 1250 + 36P_3 + 0.0435P_3^2 + 0.6H_3 + 0.027H_3^2 + 0.011P_3H_3$$

$$\text{cost}_4(H_4) = 23.4H_4$$

The domains of six decision variables are stated:
$P_1 \in [0, 150]$, $P_2 \in [81, 274]$, $P_3 \in [40, 125.8]$, $H_2 \in [0, 180]$, $H_3 \in [0, 135.6]$, $H_4 \in [0, 2695.2]$

Subjected to the power and heat balance equation constraints:

h1: $P_1 + P_2 + P_3 - P_d = 0$

h2: $H_2 + H_3 + H_4 - H_d = 0$

Dual dependency constraints of co-generation unit 2

g1: $P_2 + 0.177778H_2 - 247 \le 0$

g2: $\begin{cases} 98.8 - P_2 - 0.16985H_2 \le 0 & H_2 \in [0, 104.8] \\ -P_2 + 1.781915H_2 - 105.74468 \le 0 & H_2 \in (104.8, 180] \end{cases}$

Dual dependency constraints of co-generation unit 3

g3: $\begin{cases} P_3 - 125.8 \le 0 & H_3 \in [0, 32.4] \\ P_3 + 0.151163H_3 - 130.6977 \le 0 & H_3 \in (32.4, 135.6] \end{cases}$

g4: $\begin{cases} 44 - P_3 \le 0 & H_3 \in [0, 15.9] \\ -P_3 - 0.0067682H_3 + 45.076142 \le 0 & H_3 \in (15.9, 75] \\ -P_3 + 1.1584H_3 - 46.8812 \le 0 & H_3 \in (75, 135.6] \end{cases}$

To demonstrate the tweaking impact of the key parameters $N_L$ and $N_a$ of the TSCO, 20 trails has been executed using different knowledge library sizes and the learning agent sizes in test case, and the obtained solutions are presented in Table 1.

Table 1. Impact of parameters $N_L$ and $N_a$ on test case 1

| $N_L$ | $N_a$ | execution time (s) | Objective function value | | |
|---|---|---|---|---|---|
| | | | minimum | maximum | mean |
| 100 | 20 | 0.453 | 9259.08 | 9276.23 | 9263.34 |
| | 50 | 0.514 | 9257.08 | 9269.35 | 9260.19 |
| | 70 | 0.613 | 9257.07 | 9266.39 | 9259.82 |
| | 20 | 0.472 | 9257.10 | 9276.23 | 9261.25 |
| **150** | **50** | **0.535** | **9257.07** | **9261. 52** | **9258.11** |
| | 70 | 0.622 | 9257.07 | 9261.32 | 9258.09 |
| | 20 | 0.492 | 9257.10 | 9270.58 | 9261.18 |
| 200 | 50 | 0.542 | 9257.07 | 9262.71 | 9258.43 |
| | 70 | 0.636 | 9257.07 | 9261.21 | 9258.27 |

From Table 1, we can see that the learning agent size plays a more important role than the knowledge library size in the execution process, and the time costs mainly depend on the learning agent size. The reason for this phenomenon is that the overall calculation time of TSCO is the sum of $N_L$ and $N_a * T$. And furthermore, Table 1 suggests that global solutions can be achieve very well with appropriate time in the test case 1 when the knowledge library size and learning agent size are respectively taken as 150 and 50. In addition, there is no great significant improvement if the sizes are increased beyond those values, but it increases the execution time which is not desirable in real time problems. Therefore, the parameters $N_L$ and $N_a$ are chosen as 150 and 50 to achieve optimal performance in this test case.

The optimal solution of the test case 1 attained by the TSCO algorithm is $ 9257.07, which is showed in Table 2 with focus on the comparison of TSCO in the field of minimum fuel cost and computational time with the earlier literature algorithms.

From Table 2, the TSCO algorithm obtains global optimal fuel cost with less computational time than the traditional SCO and other methods. The reason for this phenomenon is that human has the ability of observational learning and tournament learning and has more intelligence than other smarms. Meanwhile, the Tent map local searching improves the ergodicity of the SCO. Therefore, the conclusion can be drawn that TSCO algorithm is an effective way to solving the CHPED problem.

**Table 2** Optimal solutions of CHPED for test case 1.

| Output | P1 | P2 | P3 | H2 | H3 | H4 | Min | Time(s) |
|---|---|---|---|---|---|---|---|---|
| ACSA [12] | 0.08 | 150.93 | 49 | 48.84 | 65.79 | 0.37 | 9452.2 | - |
| GA [29] | 0 | 159.23 | 40.77 | 39.94 | 75.06 | 0 | 9267.5 | - |
| RGA [30] | 0 | 158.18 | 41.82 | 37 | 78 | 0 | 9263.28 | 2.08 |
| EP [31] | 0 | 160 | 40 | 40 | 75 | 0 | 9257.1 | - |
| FA [17] | 0.0014 | 159.9986 | 40 | 40 | 75 | 0 | 9257.1 | - |
| IWO [21] | 0.0002 | 159.9998 | 40 | 40 | 75 | 0 | 9257.08 | - |
| CPSO [32] | 0 | 160 | 40 | 40 | 75 | 0 | 9257.08 | 1.42 |
| RCGA-IMM [8] | 0 | 160 | 40 | 40 | 75 | 0 | 9257.075 | - |
| HS [33] | 0 | 160 | 40 | 40 | 75 | 0 | 9257.07 | - |
| IGA-MU [34] | 0 | 160 | 40 | 39.99 | 75 | 0 | 9257.07 | - |
| SARGA [35] | 0 | 159.99 | 40.01 | 39.99 | 75 | 0 | 9257.07 | 3.76 |
| EMA [36] | 0 | 160 | 40 | 40 | 75 | 0 | 9257.07 | 0.9846 |
| TVAC_PSO [32] | 0 | 160 | 40 | 40 | 75 | 0 | 9257.07 | 1.33 |
| Direct method [37] | 0 | 160 | 40 | 40 | 75 | 0 | 9257.07 | - |
| GWO[3] | 0 | 160 | 40 | 40 | 75 | 0 | 9257.07 | 1.30 |
| MCSA[14] | 0 | 160 | 40 | 40 | 75 | 0 | 9257.07 | 1.35 |

| | | | | | | | |
|---|---|---|---|---|---|---|---|
| CSA[14] | 0 | 160 | 40 | 40 | 75 | 0 | 9257.07 | 0.59 |
| CSO[38] | 0 | 160 | 40 | 40 | 75 | 0 | 9257.07 | 1.18 |
| SCO | 0 | 160 | 40 | 40 | 75 | 0 | 9257.07 | 0.673 |
| **TSCO** | **0** | **160** | **40** | **40** | **75** | **0** | **9257.07** | **0.535** |

Considering that GA is the most classic and representative heuristic algorithm, GWO is an up-to-date bio-inspired optimization algorithm, SCO is the latest human learning-based swarm optimization algorithm, these typical heuristic algorithms GA, GWO and SCO are chosen to evaluate the performances of the TSCO for solving CHPED issues. For ease of analysis, the population sizes are set to 50 in all these algorithms. In the GA, the string length, crossover probability and mutation probability are respectively set to 72, 0.90, and 0.04.

The convergence speeds of the four algorithms in CHPED problem are shown in Fig. 2.

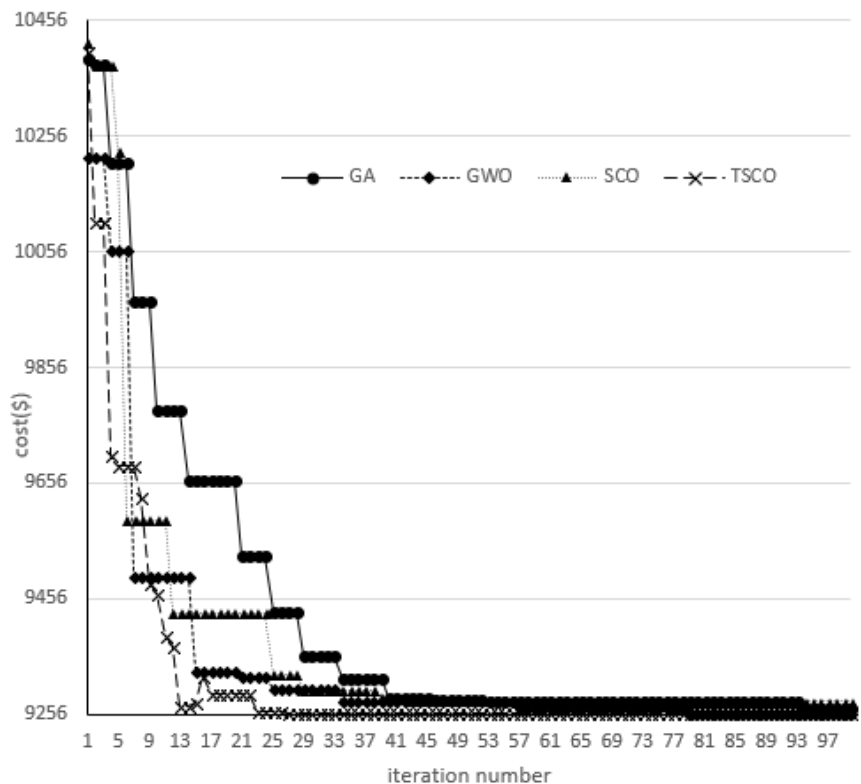

FIG. 2. Convergence comparison in case 1 by GA, GWO, SCO and TSCO

From Fig. 2, it can be observed that the TSCO has the faster convergence speed than the other three classical algorithms, especially the GA. Meanwhile, the GA doesn't achieve the optimal solution in the 100 iterations, the solution quality of the TSCO is superior to the other three algorithms. During the iteration process, because the basic constraint handling rule is that any feasible point is preferred over any unfeasible point and the point having less constraint violation is preferred among the unfeasible points, some point's objective values in the later iteration are greater than the value in the former iteration. The TSCO shows faster convergence speed and upper steady ability in the iterations of evolution than traditional SCO, GA and GWO. The reason for this phenomenon is that TSCO has the high intelligence and more global convergence performance and integrates Tent map into TSCO initialization and neighborhood searching to improve uniform distribution and increase the ergodicity of the SCO. Therefore, the conclusion can be drawn that the TSCO has good global convergence ability and Tent map is an effective way in initialization and neighborhood searching.

To reasonably compare the solution performances and stabilities of the four algorithms, 20 independent trails have been executed by using each algorithm. Fig. 3 is the boxplot of the cost difference values between the optimal value

9257.07 and the cost solution obtained by the four algorithms GA, GWO, SCO and TSCO. The vertical coordinate represents the difference values and the horizontal coordinate represents the corresponding algorithm of the boxes.

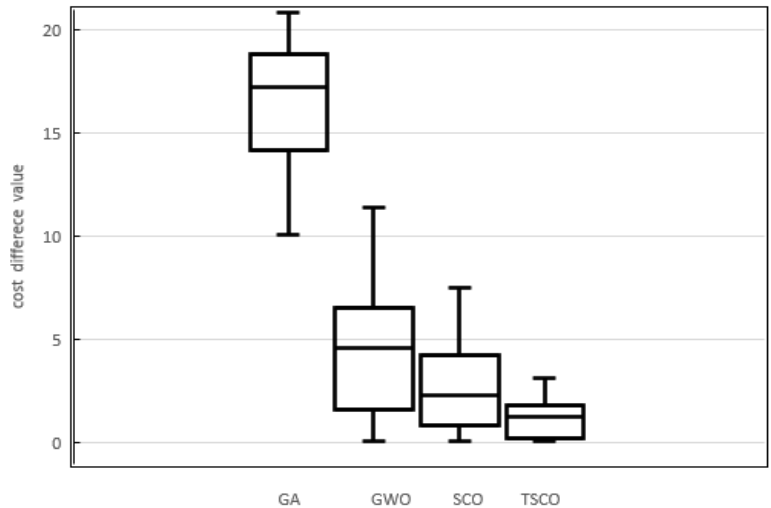

FIG. 3. The cost boxplot of four algorithms for testcase 1

It can be seen from Fig. 3 that the minimum, maximum, median, lower quartile and upper quartile of the cost difference values of the four classical algorithms are significantly different with each other. From the distribution of cost difference values, we can see that the box of TSCO is lower and shorter than those of the GA, GWO and SCO. Hence the TSCO algorithm has higher stability and better performance than the other three algorithms.

Fig.4 illustrates the constraint violation value convergence comparison of the BCH and ACR strategy in TSCO algorithm for the CHPED problem. The vertical coordinate represents constraint violation value logarithm to the base of ten, and the horizontal coordinate represents the iteration number. Because the zero logarithm does not exist, the curve of the ACR vanishes from the forty-eighth iteration. From the convergence curve, the constraint violation values of the BCH are better than those of the ACR in the first 16 iterations, because the ACR strategy introduces the relaxed quasi-feasible space. But the constraint violation values of ACR are better than those of the BCH and quickly move to zero at about forty eighth iteration. The constraint violation value of the BCH never equals zero, so the algorithm with the BCH strategy does not obtain the feasible solution within 100 iterations. The ACR strategy only needs 48 iterations from the unfeasible solutions to feasible solution area. So, the conclusion can be drawn that adaptive constraints relaxing (ACR) strategy is a very effective method to handle the equation constraints in the CHPED problem.

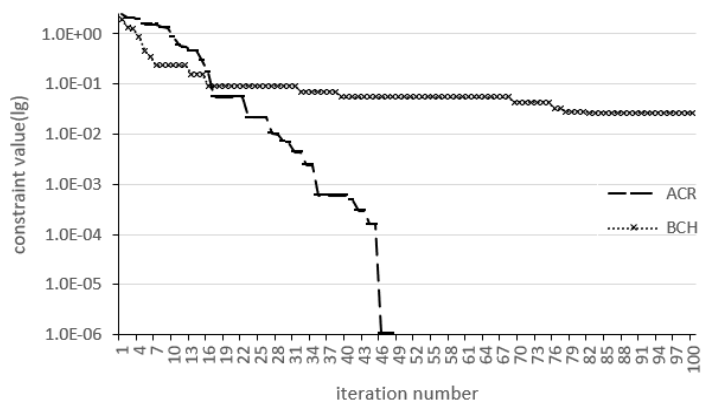

FIG. 4. Constraint violation value comparison in case 1 by ACR and BCH

## 4.2 Test case 2

In the test case 2, the valve-point effects of the power unit and power net transmission losses in the CHPED problem are considered as the same as in all the other referred literatures [3]. The test case 2 is a classical representation of sophisticated CHPDE problems and used to evaluate the algorithms in most literature. Test case 2 consists of four conventional power units, two co-generation units, and one heat unit. The test case total heat demand is 150MWth and the power demand is 600MW. The nine decision variables of the test case can be represented as $(P_1, P_2, P_3, P_4, P_5, P_6, H_5, H_6, H_7)$

The fuel cost formulation of four conventional power units, two co-generation units, and heat unit are given:

$$\text{cost}_1(P_1) = 0.008P_1^2 + 2P_1 + 25 + |100\sin(0.042 \times (10-P_1))|$$

$$\text{cost}_2(P_2) = 0.003P_2^2 + 1.8P_2 + 60 + |140\sin(0.04 \times (20-P_2))|$$

$$\text{cost}_3(P_3) = 0.0012P_3^2 + 2.1P_3 + 100 + |160\sin(0.038 \times (30-P_3))|$$

$$\text{cost}_4(P_4) = 0.001P_4^2 + 2P_4 + 120 + |180\sin(0.037 \times (40-P_4))|$$

$$\text{cost}_5(P_5, H_5) = 2650 + 14.5P_5 + 0.0345P_5^2 + 4.2H_5 + 0.03H_5^2 + 0.031P_5H_5$$

$$\text{cost}_6(P_6, H_6) = 1250 + 36P_6 + 0.0435P_6^2 + 0.6H_6 + 0.027H_6^2 + 0.011P_6H_6$$

$$\text{cost}_7(H_7) = 950 + 2.0109H_7 + 0.038H_7^2$$

$$\begin{aligned}
P_{Loss} = & (49P_1 + 14P_2 + 15P_3 + 15P_4 + 20P_5 + 25P_6)P_1 \times 10^{-7} \\
& + (14P_1 + 45P_2 + 16P_3 + 20P_4 + 18P_5 + 19P_6)P_2 \times 10^{-7} \\
& + (15P_1 + 16P_2 + 39P_3 + 10P_4 + 12P_5 + 15P_6)P_3 \times 10^{-7} \\
& + (15P_1 + 20P_2 + 10P_3 + 40P_4 + 14P_5 + 11P_6)P_4 \times 10^{-7} \\
& + (35P_5^2 + 34P_5P_6 + 39P_6^2) \times 10^{-7} \\
& + (-0.3908P_1 - 0.1297P_2 + 0.7074P_3 + 0.0591P_4 + 0.2161P_5 - 0.6635P_6) \times 10^{-3} + 0.056
\end{aligned}$$

The domains of nine decision variables are listed:
$P_1 \in [10, 75]$, $P_2 \in [20, 125]$, $P_3 \in [30, 175]$, $P_4 \in [40, 250]$, $P_5 \in [81, 247]$, $H_5 \in [0, 180]$, $P_6 \in [40, 125.8]$, $H_6 \in [0, 135.6]$, $H_7 \in [0, 60]$

Subjected to the power and heat balance equation constraints:

h1: $\sum_{i=1}^{6} P_i - P_d - P_{Loss} = 0$

h2: $\sum_{j=5}^{7} H_j - H_d = 0$

Dual dependency constraints of co-generation unit 5

g1: $P_5 + 0.177778H_5 - 247 \leq 0$

g2: $\begin{cases} 98.8 - P_5 - 0.16985H_5 \leq 0 & H_5 \in [0, 104.8] \\ -P_5 + 1.781915H_5 - 105.74468 \leq 0 & H_5 \in (104.8, 180] \end{cases}$

Dual dependency constraints of co-generation unit 6

$$g3: \begin{cases} P_6 - 125.8 \leq 0 & H_6 \in [0, 32.4] \\ P_6 + 0.151163 H_6 - 130.6977 \leq 0 & H_6 \in (32.4, 135.6] \end{cases}$$

$$g4: \begin{cases} 44 - P_6 \leq 0 & H_6 \in [0, 15.9] \\ -P_6 - 0.0067682 H_6 + 45.076142 \leq 0 & H_6 \in (15.9, 75] \\ -P_6 + 1.1584 H_6 - 46.8812 \leq 0 & H_6 \in (75, 135.6] \end{cases}$$

The simulation result obtained by the TSCO algorithm is $ 1094.2351, which is showed in Table 3 with focus on the comparison of TSCO in the field of minimum fuel cost and computational time with existing state-of-the-art algorithms.

**Table 3** Optimal solutions of CHPED for test case 2.

| Method | P1 | P2 | P3 | P4 | P5 | P6 | H5 | H6 | H7 | Min ($) | Time(s) |
|---|---|---|---|---|---|---|---|---|---|---|---|
| RCGA[8] | 74.6834 | 97.9578 | 167.2308 | 124.9079 | 98.8008 | 44.0001 | 58.0965 | 32.4116 | 59.4919 | 10667 | 6.4723 |
| PSO[32] | 18.4626 | 124.2602 | 112.7794 | 209.8158 | 98.8140 | 44.0107 | 57.9236 | 32.7603 | 59.3161 | 10613 | 5.3844 |
| EP[31] | 61.3610 | 95.1205 | 99.9427 | 208.7319 | 98.8 | 44 | 18.0713 | 77.5548 | 54.3739 | 10390 | 5.275 |
| AIS[15] | 50.1325 | 95.5552 | 110.7515 | 208.7688 | 98.8 | 44 | 19.4242 | 77.0777 | 53.4981 | 10355 | 5.2956 |
| CPSO [32] | 75 | 112.38 | 30 | 250 | 93.2701 | 40.1585 | 32.5655 | 72.6738 | 44.7606 | 10325.399 | 3.29 |
| DE [19] | 44.2118 | 98.5383 | 112.6913 | 209.7741 | 98.8217 | 44 | 12.5379 | 78.3481 | 59.1139 | 10317 | 5.26 |
| BCO [13] | 43.9457 | 98.5888 | 112.932 | 209.7719 | 98.8 | 44 | 12.0974 | 78.0236 | 59.879 | 10317 | 5.1563 |
| ECSA [39] | 53.7610 | 98.5039 | 112.5996 | 209.7993 | 93.0872 | 40.2022 | 33.6571 | 72.6890 | 43.6539 | 10121.9466 | 0.434 |
| KHA [40] | 46.3835 | 104.1223 | 64.3729 | 246.1853 | 98.9736 | 40.7401 | 0 | 66.71 | 83.29 | 10111.1501 | 2.2094 |
| EMA[36] | 52.6847 | 98.5398 | 112.6734 | 208.8158 | 93.8341 | 40 | 29.242 | 75 | 45.7579 | 10111.0732 | 2.0654 |
| TLBO[23] | 45.266 | 98.5479 | 112.6786 | 209.8284 | 94.4121 | 40.0062 | 25.8365 | 74.9970 | 49.1666 | 10094.8384 | 2.86 |
| CSO[38] | 45.4909 | 98.5398 | 112.6734 | 209.8158 | 94.1838 | 40 | 27.1786 | 75 | 47.8214 | 10094.1267 | 3.09 |
| SCO | 58.7268 | 98.5398 | 112.6735 | 209.8158 | 81 | 40 | 92.0061 | 45.0590 | 12.9349 | 10226.8556 | 0.31 |
| **TSCO** | **45.5231** | **98.5385** | **112.5798** | **209.8159** | **94.2261** | **40** | **28.6947** | **74.9981** | **46.3072** | **10094.2351** | **0.29** |

From Table 3, it can be observed that the CPU time of TSCO algorithm is much less than all other compared algorithms. In addition, the TSCO algorithm obtains less fuel cost than the algorithm reported in the existing literature [6, 11, 13, 17, 21, 29, 30, 34, 36, 37, 38], and almost equal fuel cost with the CSO [36] algorithm. As shown in Table 2, the power production and the heat production are 600.6834 MW and 150 WMth, respectively. The power net transmission loss is 0.6834 MW. Transparently, the output result completely fulfills the heat and power demands. The reason for this phenomenon is that human has the ability of observational learning and tournament learning and has more intelligence than other smarms. Meanwhile, the Tent map local searching improves the ergodicity of the SCO. Therefore, the result of the CPU computation time of TSCO is much less than that of the classical SCO. Therefore, the conclusion can be drawn that TSCO algorithm is an effective way to solving CHPED problem.

The cost convergence quality acquired by TSCO and SCO is presented in Fig.5.

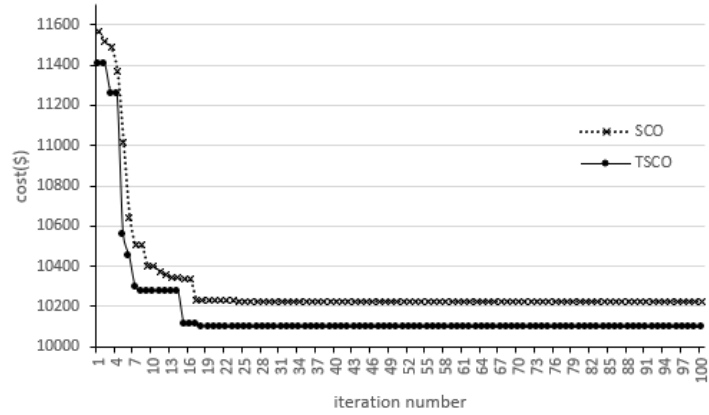

FIG. 5. Convergence comparison in case 2 by SCO and TSCO

As shown in Fig. 5, the TSCO shows faster convergence speed and upper steady ability than traditional SCO in the iterations of evolution. The reason for this phenomenon is that TSCO algorithm introduces Tent map to improve uniform distribution and avoid the prematurity convergence of the algorithm. Therefore, the conclusion can be drawn that Tent map is very effective in initialization and neighborhood local searching.

The constraint violation value convergence comparison of the BCH and ACR strategy of the TSCO algorithm in the sophisticated CHPED problem is illustrated in Fig. 6. From the convergence curve in Fig. 6, the constraint violation values of the BCH are better than those of the ACR in the first 8 iterations. The reason for this phenomenon is that the ACR strategy introduces the relaxed quasi-feasible space and obtains the unfeasible solution in the early iterations period. But the constraint violation values of ACR are better than those of the BCH in the later iterations period and quickly move to zero after about forty seventh iteration. The constraint violation value of the BCH never equals zero, so the algorithm with the BCH strategy does not obtain the feasible solution within 100 iterations. The ACR strategy only needs 47 iterations from the unfeasible solutions to feasible solution area. The reason for this phenomenon is that the ACR strategy enhances the movement to the feasible space in the later iterations period. So, the conclusion can be drawn that adaptive constraints relaxing (ACR) strategy is very effective to handle the equation constraints in the sophisticated CHPED problem.

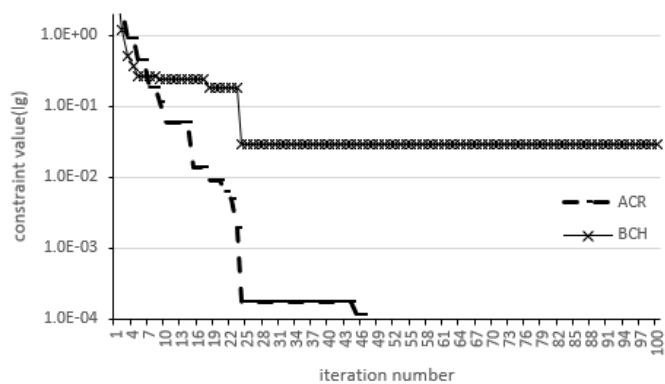

FIG. 6. Constraint violation value comparison in case 2 by ACR and BCH

## 5. Conclusions

To acquire the best utilization of CHP units, a social cognitive optimization algorithm with tent map for the CHPED problem is presented. The presented algorithm is verified by two classical test cases. The conclusions can be drawn as follows:

(1) TSCO algorithm can effectively settle the CHPED problem, and it has better solution characteristic, less CPU computation time and upper convergence speed than many intelligence optimization algorithms.

(2) Adaptive constraints relaxing rule is adopted in handling the equation constraint processing in heat and power balance constraints, which is very effective in the CHPED problem.

(3) Integrating Tent map into SCO initialization and neighborhood searching obtains better optimization than the classical stochastic method in the SCO algorithm. Furthermore, the method of the chaotic search may be applied to any similar NLP problems.

Future work will focus on integrating renewable energy resources, such as wind and solar generation, into CHP systems to construct an integrated energy system. Besides, it is interesting to investigate more realistic modeling techniques such as load and renewable generation uncertainties and energy storage units [41] to improve the practicality of our approach.

**Acknowledgments**

The work is supported by the special fund for key discipline construction of general institutions of higher learning from Shaanxi Province. We thank Dr. Xie for sharing the source code of social cognitive optimization algorithm on the website: http://www.wiomax.com/team/xie/social-cognitive-optimization-sco-project-portal/.